\documentclass[review]{elsarticle}

\usepackage{lineno}
\usepackage{amssymb}
\modulolinenumbers[5]

\newtheorem{thm}{Theorem}
\newproof{pf}{Proof}

\begin{document}

\begin{frontmatter}

\title{Stratified Monte Carlo simulation of Markov chains}

\author[mymainaddress]{Rana Fakhereddine}
\ead{	ranafakhreddine@hotmail.com}

\author[mymainaddress]{Rami El Haddad}
\ead{rami.haddad@usj.edu.lb}

\author[mysecondaryaddress]{Christian L\'ecot\corref{mycorrespondingauthor}}
\cortext[mycorrespondingauthor]{Corresponding author}
\ead{Christian.Lecot@univ-savoie.fr}

\address[mymainaddress]{Universit\'e Saint-Joseph, Facult\'e des Sciences, BP 11-514 Riad El Solh, Beyrouth 1107 2050, Liban}
\address[mysecondaryaddress]{Universit\'e de Savoie, LAMA, UMR 5127, Campus universitaire,
73376 Le Bourget-du-Lac, France}

\begin{abstract}
We present several Monte Carlo strategies for simulating discrete-time Markov chains with continuous multi-dimensional state space; we focus on stratified techniques. We first analyze the variance of the calculation of 
the measure of a domain included in the unit hypercube, when stratified samples are used. We then show that each step of the simulation of a Markov chain can be reduced to the numerical integration of the indicator function of a subdomain of the unit hypercube. Our approach for Markov chains simulates $N$ copies of the chain in parallel using stratified sampling and the copies are sorted after each step, according to their successive coordinates. We analyze variance reduction on examples of pricing of European and Asian options: enhanced efficiency of stratified strategies is shown.
\end{abstract}


\end{frontmatter}

\linenumbers

\section{Introduction}
\label{sec:1}

Many real-life systems can be modeled using Markov chains. Fields of application are
queueing theory, telecommunications, option pricing, etc. In most interesting
situations, analytic formulas are not available and the state space
of the chain is so large that classical numerical methods would require
a considerable computational time and huge memory capacity. So
Monte Carlo (MC) simulation becomes the standard way of estimating
performance measures for these systems. A drawback of MC methods is
their slow convergence, with respect to the number of random points used. 
Various techniques have been developed, in order to 
reduce the variance of the approximation, including stratified sampling and Latin hypercube 
sampling \cite{Fis96:MonteCarlo,ES00:Approximating,Gla04:MonteCarlo}.

It is shown in a series of papers \cite{LT04:Quasi,LT04:Comparison,ELL08:Quasi,
ELLN10:Quasi} that each step of a MC simulation of a Markov chain amounts to approximating 
the measure of a subdomain of the $s$-dimensional unit hypercube $I^s := [0,1)^s$. 
The techniques presented here use stratified samples for calculating this approximation. 

Among stratification strategies, we first consider the simple approach (SMC): the unit hypercube is divided into $N$ subcubes having the same measure, and one random point 
is chosen in each subcube. For Latin hypercube sampling (LHS), the projections of the points on each coordinate axis are evenly distributed: one projection in each of the $N$ subintervals that uniformly 
divide the unit interval $I $. Then we propose an hybrid method between SMC and LHS, that has properties of both approaches, with one random point in each subcube and one projection in each 
subinterval; we call this technique \emph{Sudoku} Sampling (SS) due to the properties of the points recalling a Sudoku grid.

The improved accuracy of stratified methods may be lost for problems in which we have to approximate 
the measure of subdomains with irregular boundaries. It is necessary to take special measures to make optimal use of the greater uniformity associated with stratified samples. This is achieved in \cite{LLT06:Randomized,LLT08:Randomized,LLL09:Array} through the additional effort of reordering the copies of the chain at each time step. This type of sorting was initiated by \cite{Lec89:Direct} in the context of quasi-Monte Carlo (QMC) methods.

This paper is organized as  follows. In Section \ref{sec:Integration} we present SMC, LHS and SS methods for numerical integration. We recall variance bounds for SMC and SS and we establish a new bound for the variance of LHS approach in the restrictive case of the approximate calculation of the measure of an interval in dimension $s$. In Section \ref{sec:Simulation}, we propose a MC simulation of Markov chains using stratified samples in the context of discrete Markov chains with continuous multi-dimensional state space. The results of  numerical experiments are presented in Section \ref{sec:Illustrations}. We compute the values of European and Asian options and we compare the variance of the results and the efficiency of the approaches. It is shown that both SMC and SS strategies outperform MC or LHS approaches. Finally, we give some perspectives for future work.

\section{Numerical integration}
\label{sec:Integration}

Let $s \geq 1$ be a given dimension; then $I^s$ is the $s$-dimensional half-open unit hypercube and $\lambda_s$ denotes the $s$-dimensional Lebesgue measure. If $g$ is a square-integrable function defined on $I^s$, we want to approximate
\begin{equation}
\label{eq:int}
\mathcal{I} := \int_{I^s} g(x) d\lambda_s(x). 
\end{equation}
For the usual MC approximation, $\{U_1, \ldots , U_N \}$ are independent random variables uniformly distributed over $I^s$. Then 
\begin{equation}
\label{eq:intMC}
X := \frac{1}{N} \sum_k g (U_k)
\end{equation}
is an unbiased estimator of $\mathcal{I}$. When $g=1_A$, for some measurable $A \subset I^s$, one has 
\begin{equation}
\mathrm{Var}(X) = \frac{1}{N} \lambda_s( A ) ( 1 - \lambda_s( A ))  \leq \frac{1}{4N}.
\end{equation}
A simple stratified sampling (SMC) method was proposed in \cite{Hab66:Modified} and further analyzed in \cite{CD89:Problem}. For $N = n^s$, put
\begin{equation}
\label{eq:Jell}
J_{\ell} := \prod_{i=1}^s \left[ \frac{\ell_i-1}{n}, \frac{\ell_i}{n} \right), \quad 1 \leq \ell_1 \leq n, 
\ldots , 1 \leq \ell_s \leq n.
\end{equation}
Let $\{V_{\ell} : 1 \leq \ell_1 \leq n, 
\ldots , 1 \leq \ell_s \leq n \}$ be independent random variables,  
with $V_{\ell}$ uniformly distributed over $J_{\ell}$. Then
\begin{equation}
\label{eq:EstSMC}
Y := \frac{1}{N} \sum_{\ell} g (V_{\ell})
\end{equation}
is another unbiased estimator of $\mathcal{I}$. In \cite{EFL13:Stratified}, we have analysed the following case:  we consider a function $f : \overline{I}^{s-1} \to \overline{I}$ and we define 
\begin{equation}
\label{eq:Af}
A_f := \{ (u^{\prime},u_s) \in I^s : u_s < f(u^{\prime}) \}.
\end{equation}
Then for $g =1_{A_f}$ we obtain
\begin{equation}
\label{eq:VarSMC}
\mathrm{Var}(Y) \leq \left( \frac{s-1}{4} V(f) + \frac{1}{2} \right) \frac{1}{N^{1+1/s}},
\end{equation}
if $f$ is of bounded variation $V(f)$ in the sense of Hardy and Krause (we refer to \cite{Nie92:Random} for this concept). 

Latin hypercube sampling (LHS) was introduced in \cite{MBC79:Comparison} and examined studiously in \cite{Ste87:Large,Owe97:MonteCarlo}. Let 
\begin{equation}
\label{eq:IntIk}
I_k := \left[ \frac{k-1}{N},\frac{k}{N} \right), \quad 1 \leq k \leq N
\end{equation}
and $\{V_1^i, \ldots , V_N^i \}$ be independent random variables,  where $V_k^i$ is uniformly distributed over $I_k$. If $\{ \pi^1, \ldots , \pi^s \}$ are independent random permutations of $\{1, \ldots , N \}$ and
$W_k := ( V_{\pi^1(k)}^1, \ldots , V_{\pi^s(k)}^s )$, then each $W_k$ is uniformly distributed over $I^s$. 
Consequently,
\begin{equation}
\label{eq:intLHS}
Z := \frac{1}{N} \sum_k g (W_{k})
\end{equation}
is another unbiased estimator of $\mathcal{I}$. 

We have proposed in \cite{EFLV14:Extended} a combination of SMC and LHS: we construct $N = n^s$ random points in $I^s$ such that in every interval
$I^{i-1} \times I_k \times I^{s-i}$ (for $1 \leq i \leq s$ and $1 \leq k \leq N$) or $J_{\ell}$ 
(for $1 \leq \ell_1 \leq n, \ldots , 1 \leq \ell_s \leq n$) lies only one point of the set (property $\mathcal{P}$). This is achieved as follows. If $x := (x_1, \ldots , x_s)$, we put 
$\hat{x}_i := (x_1, \ldots , x_{i-1},x_{i+1}, \ldots , x_s)$.
Let $\sigma^1, \ldots , \sigma^s$
be random bijections $\{1, \ldots , n\}^{s-1} \to \{ 1, \ldots , n^{s-1} \}$
and $\{ U_{\ell}^i : 1 \leq i \leq s, 1 \leq \ell_1 \leq n, \ldots , 1 \leq \ell_s \leq n \}$ be random variables uniformly distributed on $I$; all these variables are assumed to be mutually independent. We put
\begin{equation}
W_{\ell}^{\ast} = \Big( \frac{\ell_1-1}{n} + \frac{\sigma^1(\hat{ \ell}_1)-1+U^1_{\ell}}{N}, \ldots , \frac{\ell_s-1}{n} + \frac{\sigma^s(\hat{ \ell}_s)-1+U^s_{\ell}}{N} \Big).
\label{eq:defWelli}
\end{equation}
The point set $\{W_{\ell}^{\ast} : 1 \leq \ell_1 \leq n, \ldots , 1 \leq \ell_s \leq n \}$ has property $\mathcal{P}$. 
If $Z^{\ast}$ is defined by 
\begin{equation}
\label{eq:intSS}
Z^{\ast} := \frac{1}{N} \sum_{\ell} g (W_{\ell}^{\ast}),
\end{equation}
it is an unbiased estimator of $\mathcal{I}$. The following variance bound is established in \cite{EFLV14:Extended}. Let $A \subset I^s$ be such that, for all $i$, with $1 \leq i \leq s$,
\begin{equation}
\label{eq:defA}
A = \{ (u_1, \ldots , u_s) \in I^s : u_i < f_i(\hat{u}_i) \},
\end{equation}
where $f_i$ are Lipschitz continuous functions 
$\overline{I}^{s-1} \to \overline{I}$. Then, for $g = 1_A$, we have
\begin{equation}
\label{eq:VarSS}
\mathrm{Var}(Z^{\ast}) \leq \Big( \frac{\kappa+2}{4}+2s(\kappa+2)^2 \Big) \frac{1}{N^{1+1/s}}, 
\end{equation}
where $\kappa$ is a Lipschitz constant (for the maximum norm) for all the $f_i$.
We prove a similar result for LHS, in a very restrictive case.

\begin{thm}
If $A$ is a subinterval of $I^s$, then the variance of the LHS approximation
\[
Z=\frac{1}{N} \sum_k  1_A(W_k)
\]
satisfies (for $N \geq 3$)
\[
\mathrm{Var}(Z) \leq \frac{1}{N} \lambda_s(A) ( 1 - \lambda_s(A)).
\]
\end{thm}

\begin{pf}
We have
\begin{eqnarray*}
\mathrm{Var}(Z) & = & \frac{1}{N^2} \sum_k \mathrm{Var}(1_A(W_k))
+ \frac{1}{N^2} \sum_{ k \neq k^\prime} \mathrm{cov}(1_A(W_k),1_A(W_{k\prime})) \\
& = & \frac{1}{N} \lambda_s(A)( 1 - \lambda_s(A))
+ \frac{1}{N^2} \sum_{ k \neq k^\prime} \mathrm{cov}(1_A( W_k),1_A(W_{k^\prime})).
\end{eqnarray*}
For $\overline{j} := (j_1, \ldots , j_s)$ with $1 \leq j_1 \leq N, \ldots , 1 \leq j_s \leq N$, we set
$\overline{J}_{\overline{j}} := \prod_{i=1}^s I_{j_i}$.  If $k \neq k^{\prime}$, then 
\begin{eqnarray*}
\mathrm{cov}(1_A(W_k),1_A(W_{k^\prime}))  & = &
\frac{N^s}{(N-1)^s} \sum_{j_1 \neq j_1^\prime} \cdots
\sum_{j_s \neq j_s^\prime}  \lambda_s(A \cap \overline{J}_{\overline{j}}) 
\lambda_s(A\cap \overline{J}_{\overline{j}^{\prime}})  \\
& & - (\lambda_s(A))^2.
\end{eqnarray*}
We may assume that $A$ is a closed interval:
\[
A := \prod_{i=1}^s \Big[ \frac{m_i-x_i^{-1}-1}{N} , \frac{m_i+n_i+x_i^{+1}-1}{N} \Big],
\]
with $1 \leq m_i,n_i \leq N$ and $x_i^{-1},x_i^{+1} \in I$.
Let us note $[1,s] := \{1,2, \ldots , s \}$ and, for $H \subset [1,s]$, denote $H^c := [1,s] \setminus H$. Then
\begin{eqnarray*}
\lefteqn{ \sum_{j_1 \neq j_1^\prime} \cdots
\sum_{j_s \neq j_s^\prime}  \lambda_s(A \cap \overline{J}_{\overline{j}}) 
\lambda_s(A\cap \overline{J}_{\overline{j}^{\prime}}) } \\
& = & \frac{1}{N^{2s}} \sum_{H \subset [1,s]} \sum_{\epsilon_i = \pm 1, i \in H^c} \prod_{h \in H} n_h
(n_h+x_h^{-1} + x_h^{+1} -1) \prod_{i \in H^c} x_i^{\epsilon_i}
(n_i+ x_i^{-\epsilon_i} ) \\
& = & \frac{1}{N^{2s}} \prod_{i=1}^s ( n_i (n_i + x_i^{-1} + x_i^{+1}-1) + (n_i + x_i^{-1} ) x_i^{+1}
+ (n_i + x_i^{+1} ) x_i^{-1} ).
\end{eqnarray*}
Hence
\begin{eqnarray*}
\lefteqn{ (N^2(N-1))^s \mathrm{cov}(1_A(W_k),1_A(W_{k^\prime})) } \\
&  = & \prod_{i=1}^s N( n_i (n_i + x_i^{-1} + x_i^{+1}-1) + (n_i + x_i^{-1} ) x_i^{+1}
+ (n_i + x_i^{+1} ) x_i^{-1} ) \\
& & -  \prod_{i=1}^s (N-1)( n_i + x_i^{-1} + x_i^{+1} )^2.
\end{eqnarray*}
Since
\begin{eqnarray*}
\lefteqn{ (N-1) ( n_i+ x_i^{-1} + x_i^{+1} )^2 } \\
& & - N ( n_i (n_i + x_i^{-1} + x_i^{+1}-1) + (n_i + x_i^{-1} ) x_i^{+1} + (n_i + x_i^{+1} ) x_i^{-1} )  \\
& = & \frac{1}{2} (N-2) \Big( x_i^{-1} + x_i^{+1} - \frac{2n_i}{N-2} \Big)^2 
+ \frac{N}{2} (x_i^{-1} - x_i^{+1} )^2 + \frac{n_iN(N-2-n_i)}{N-2}, 
\end{eqnarray*}
we obtain $\mathrm{cov}(1_A(W_k),1_A(W_{k^\prime}))  \leq 0$ and the result follows.
\end{pf}

\section{Simulation of Markov chains}
\label{sec:Simulation}

In this section, we use the previous stratification techniques for Markov chains simulation. 

\subsection{Markov chain setting and Monte Carlo simulation}

Let $s \in \mathbb{N}^{\ast}$; we consider an homogeneous Markov chain $\{ X_p,\, p \in \mathbb{N}\}$
with state space $E\subset\mathbb{R}^s$,
evolving according to the stochastic recurrence: for $p \geq 0$
\begin{equation}
\label{eq:chaintransition}
X_{p+1} = \varphi_{p+1}(X_p,U_{p+1}).
\end{equation}
Here $\{ U_p,\, p \geq 1 \}$ is a sequence of i.i.d. uniform random variables over $I^d$
(for $d \in\mathbb{N}^*$) and each $\varphi_{p+1} : E \times I^d \to E$
is a measurable map.
The distribution $P_0$ of $X_0$ is known, and our aim is to approximate the distribution $P_p$ of $X_p$.
The standard iterative Monte Carlo scheme proceeds as follows.    
A large number $N$ of samples $x_k^0, \ 1 \leq k \leq N$ are drawn from the  initial distribution $P_0$; then we generate $N$ sample paths of the chain as follows. For $p \geq 0$ and for each $k \in \{1,2, \ldots , N\}$
\begin{equation}
\label{eq:MCtransition}
x_k^{p+1} = \varphi_{p+1}(x_k^p,u_k),
\end{equation}
where $\{ u_k, 1 \leq k \leq N \}$ are  pseudo-random numbers 
simulating i.i.d. uniform random variables over $I^d$, independent from all variables introduced previously.
QMC variants have been proposed to improve the accuracy of the method \cite{LLT06:Randomized,ELLN10:Quasi}. The pseudo-random numbers $u_k$ are replaced with quasi-random numbers; in order to benefit from the great uniformity of quasi-random points, one possibility is to sort the states $x_k^p$ by position in every step. Since QMC methods do not give confidence intervals, randomized QMC algorithms have also been introduced in \cite{LLT06:Randomized,LLT08:Randomized,LLL09:Array}, with randomized quasi-random points. In the present paper, we propose a scheme using the sampling strategies presented in section \ref{sec:Integration}. 

\subsection{Stratified algorithm}

Let $\mathcal{M}_+(E)$ denote the set of all nonnegative measurable functions on $E$. From (\ref{eq:chaintransition}), we obtain
\begin{equation}
\label{eq:caractransition}
\forall f \in \mathcal{M}_+(E) \quad \int_E f (x) d P_{p+1}(x) = \int_{I^d} \int_E f \circ \varphi_{p+1}(x,u^{\prime \prime})
d P_p(x) du^{\prime \prime}.
\end{equation}

Let $n \geq 2$ be an integer and put $N:=n^{s+d}$. For each $p \geq 0$, we are looking for an
approximation of $P_p$ of the form
\begin{equation}
\label{eq:defapprox}
\widehat{P}_p := \frac{1}{N} \sum_k \delta(x-x_k^p),
\end{equation}
where $\Xi^p := \{ x_1^p, \ldots, x_N^p \}$ is a subset of $E$ to be determined.
We first sample a point set  $\Xi^0$ of $N$ states from the initial probability distribution $P_0$. Once we have calculated a point set $\Xi^p$ such that $\widehat P_p$ \emph{approximates} $P_p$, we compute $\Xi^{p+1}$ in two steps: we first sort the states of $\Xi^p$ according to their successive coordinates, then we perform a numerical integration using a stratified sample.  

\paragraph*{Step 1: Relabeling the states} 

We label the states $x_m^p$ using a multi-dimensional index $m = (m_1, \ldots ,m_s)$ with $1\leq m_1 \leq n, \ldots , 1\leq m_{s-1} \leq n, 1\leq m_s \leq n^{1+d}$, such that:
\begin{description}
\item[] if $m_1<m^{\prime}_1$ then $x_{m,1}^p \leq x_{m^{\prime},1}^p$,
\item[] if $m_1=m^{\prime}_1,m_2<m^{\prime}_2$ then $x_{m,2}^p \leq
x_{m^{\prime},2}^p$,
\item[] $\cdots$
\item[] if $m_1=m^{\prime}_1,\ldots,m_{s-1}=m^{\prime}_{s-1}, m_s<m^{\prime}_s$
then $x_{m,s}^p \leq x_{m^{\prime},s}^p$.
\end{description}
In the case $s=1$, this reduces to simply sort the states by increasing order. 
If $s\geq 2$, the $N$ states are first sorted in $n$ batches of size $N/n$ according to their first coordinates; then each batch is sorted in subgroups of $n$  
batches of size $N/n^2$ by order of the second coordinates, and so on. 
At the last step of the sorting, subgroups of size
$n^{d+1}$ are ordered according to the last coordinate of the state.
This type of nested sorting was introduced in \cite{LC98:Quasi} for the QMC simulation of the Boltzmann equation: since the algorithm is described by a series of numerical integrations, the sorting tends to reduce the number of the jumps of the integrand.

\paragraph*{Step 2: Using stratified samples for transition} 

We define a probability measure $\widetilde{P}_{p+1}$ on $E$ by replacing $P_p$ with $\widehat{P}_p$ in eq. \ref{eq:caractransition}:
\begin{equation}
\label{eqInt}
 \int_E f(x) d \widetilde{P}_{p+1}(x) := \int_{I^d} \int_E f \circ \varphi_{p+1}(x,u^{\prime\prime}) d \widehat{P}_p(x) du^{\prime\prime}, \quad f \in \mathcal{M}_+(E).
\end{equation}
To obtain a uniform approximation of $P_{p+1}$, similar to (\ref{eq:defapprox}), we use a quadrature with stratified samples: let $\{ w_\ell : 1\leq \ell_1\leq n,\ldots, 1\leq \ell_{s+d} \leq n\}$ be pseudo-random numbers 
simulating stratified variables on $I^{s+d}$ as described in section \ref{sec:Integration}, independent from all variables introduced previously. For $m = (m_1, \ldots ,m_s)$ with $1\leq m_1 \leq n, \ldots , 1\leq m_{s-1} \leq n, 1\leq m_s \leq n^{1+d}$, let $1_m$ be the indicator function of the interval 
$\prod_{i=1}^{s-1} [(m_i-1)/n ,m_i/n) \times [(m_s-1)/n^{1+d} ,m_s/n^{1+d})$. 
For $f \in \mathcal{M}_+(E)$, denote
\begin{equation}
\label{eq:defCf}
C^p f(u) := \sum_{m}  1_m (u^{\prime}) f \circ \varphi_{p+1} (x_m^p,u^{\prime\prime}),
\quad u =(u^{\prime},u^{\prime\prime}) \in I^s \times I^d.
\end{equation}
Then we have:
\begin{equation}
\label{eq:transintegr}
\forall f \in \mathcal{M}_+(E) \quad \int_E f(x) d \widetilde{P}_{p+1}(x) = \int_{I^{s+d}} C^p f (u) du.
\end{equation}
We obtain  $\widehat P_{p+1}$ by
\begin{equation}
\label{eq:defStratapprox}
 \int_E f(x) d \widehat P_{p+1}(x) := \frac{1}{N} \sum_{\ell} 
C^p f (w_{\ell}), \quad f \in \mathcal{M}_+(E).
\end{equation}
The second step of the algorithm may be written as follows. For $u \in I^{s+d}$ let $u^{\prime} := (u_1, \ldots u_s)$ and $u^{\prime \prime} := (u_{s+1}, \ldots , u_{s+d})$; for
$u^{\prime}\in I^s$, let
$m(u^{\prime}) := (1+\lfloor n u_1 \rfloor,\ldots,1+\lfloor n u_{s-1} \rfloor,1+ \lfloor n^{1+d}u_s \rfloor)$. Then
\begin{equation}
\label{eq:Strattransition}
x_{\ell}^{p+1} = \varphi_{p+1}(x_{m(w^{\prime}_{\ell})}^p,w^{\prime \prime}_{\ell}). 
\end{equation}
(compare with eq. \ref{eq:MCtransition}). Here the states are labeled using a multi-dimensional index $\ell = (\ell_1, \ldots ,\ell_{s+d})$ with $1\leq \ell_1\leq n,\ldots, 1\leq \ell_{s+d} \leq n$. The first $s$ components of $w_{\ell}$ are used to select the state of the chain that perform a transition, while the remaining $d$ components are used to determine the new state.

\section {Numerical illustrations}
\label{sec:Illustrations}

In this section, we compare the stratified strategies with the standard MC scheme in numerical experiments

\subsection {Pricing a European call option}

In the Black-Scholes model and under the risk-neutral measure, the asset price  $S_t$ at time $t$ obeys the stochastic differential equation:
$dS_t = r S_t dt + \sigma S_t d B_t$,
where $r$ is the risk-free interest rate, $\sigma$ the
volatility parameter and $B$ is a standard Brownian motion. The solution of this equation is given by
\begin{equation}
\label{eq:St}
S_{t}=S_{0} \exp \left((r-\sigma^2/2) t +\sigma B_t \right).
\end{equation}
Let $T$ be the maturity date and $K$ the strike price. We want to estimate 
the value of the call option: $C_{\mathcal{E}} = {\mathrm{e}}^{-rT} \mathbb{E} [(S_T -K)_+]$.
To formulate the problem  as a Markov chain, we discretize the interval $[0,T]$ 
using  observation times $0=t_0<t_1< \cdots <t_P=T$. The discrete version of 
(\ref{eq:St}) can be written as: for $p \geq 0$
\begin{equation}
\label{eq:Stp}
S_{t_{p+1}}=S_{t_p} \exp \left((r-\sigma^2/2)\Delta t_{p+1} +\sigma (B_{t_{p+1}} - B_{t_p} ) \right),
\end{equation}
where $\Delta t_{p+1} := t_{p+1}-t_p$. 
 
\begin{figure} 
\begin{center}
\includegraphics[width=0.9\textwidth,height=.4\textheight]{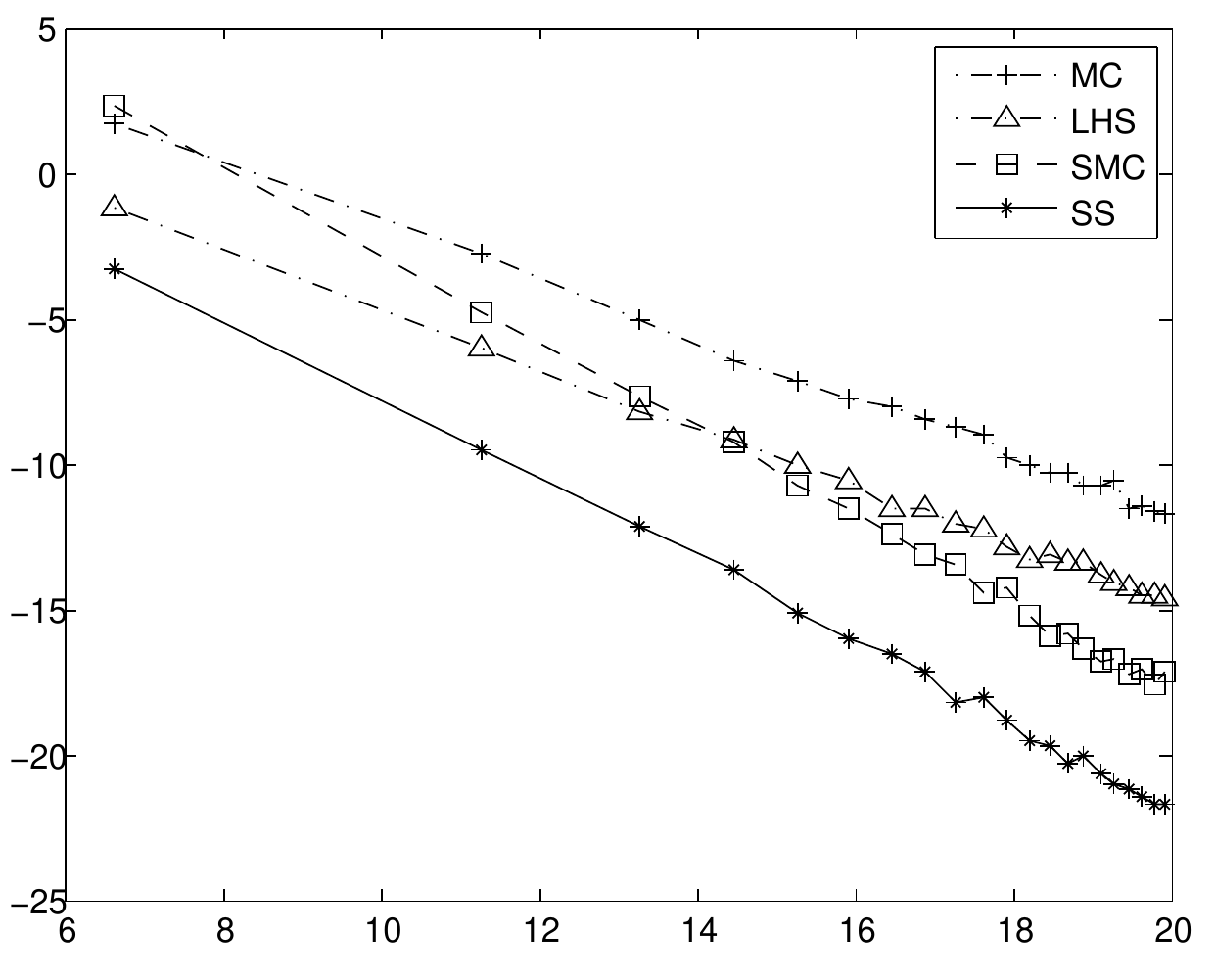} 
\end{center}
\caption{European option. Sample variance of $100$ copies of the calculation of 
$C_{\mathcal{E}}$ as a function of $N$. MC ($+$), LHS ($\triangle$), SMC ($\square$) 
and SS ($*$) outputs, in log-log scale (base $2$).}
\label{fig:europ_var}
\end{figure}

In this example $s=d=1$.  We choose the following parameters: 
$S_0=100$, $K=90$,   $r=0.06$,  $\sigma=0.2$,  $T=1$,  $P= 100$ 
 and $\Delta t_p = T/P$, for $1\leq p\leq P$.
We want to compare the variances of the MC, LHS, SMC and SS estimators of $C_{\mathcal{E}}$.  
We replicate the calculation independently $100$ times and we compute the sample variance. 
Figure \ref{fig:europ_var} shows the results as functions of $N$, for $N=10^2, 50^2,100^2,150^2,\ldots, 1000^2$, in log-log scale (base $2$). 

It is clear that SMC and SS produce smaller variances than MC and LHS (for the same $N$). When comparing the results of SMC and SS, we can see that the later approach outperforms the former. At each step of the SS algorithm, the mapping
$\ell \in \{1, \ldots , n \}^2 \to m(w_{\ell}^{\prime}) \in \{1, \ldots , n^2 \}$ is one-to-one, so that each state is considered exactly once for a transition. 

Assuming that ${\mathrm{Var}}=\mathcal{O}(N^{-\alpha})$, linear regression is used to evaluate $\alpha$. The outputs are listed in Table \ref{tab:europ}. The convergence rates are close to those established for numerical integration in dimension $2$.

\begin{table}
\begin{center}
\begin{tabular}{cccc}
\hline
MC& LHS   & SMC   & SS \\
\hline 
 1.01  &  1.01  & 1.51  & 1.42 \\
\hline
\end{tabular}
\end{center}
\caption{European option: order $\alpha$ of the variance of the calculation of $C_{\mathcal{E}}$.}
\label{tab:europ}
\end{table}   
  
\begin{figure}
\begin{center}
\includegraphics[width=0.9\textwidth,height=.4\textheight]{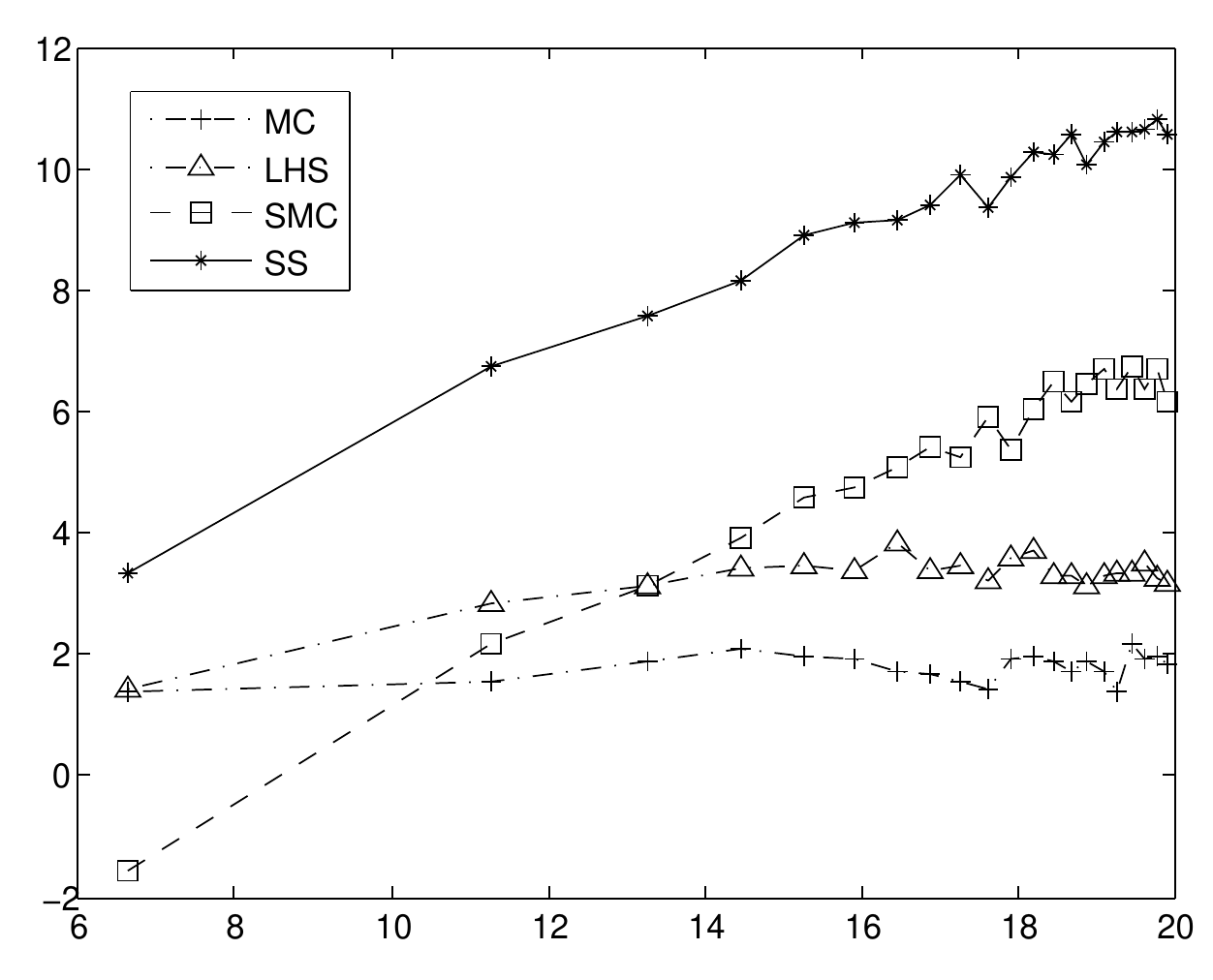} 
\end{center}
\caption{European option: efficiency of $100$ copies of the calculation of 
$C_{\mathcal{E}}$ as a function of $N$. Comparison of MC ($+$), LHS ($\triangle$), SMC ($\square$) 
and SS ($*$) outputs, in log-log scale (base $2$).}
\label{fig:europ_eff}
\end{figure}

Since we use techniques that may reduce the variance at the expense of an increase in computation time, we compare the \emph{efficiency} of the approaches. The efficiency as defined in \cite{Lec94:Efficiency} is the inverse of the product of the variance by the CPU time. It has the property that it is independent of the number $N$ of states for a naive MC estimator. The results are displayed in Figure \ref{fig:europ_eff} and show the benefits of both SMC and SS techniques.

\subsection {Pricing an Asian option}

We consider the pricing of an Asian option on a single asset. The asset price  $S_t$ at time $t$ satisfies (\ref{eq:St}) and the value of the call option with strike price $K$ at maturity date $T$ is given by 
\begin{equation}
\label{eq:CA}
C_{\mathcal{A}}={\mathrm{e}}^{-rT} \mathbb{E} \Bigl[ \Big( \Bigl(\prod_{p=1}^P S_{t_p}\Bigr)^{1/P}-K\Bigr)_+\Bigr],
\end{equation}
where   $0=t_0<t_1< \cdots <t_P=T$ are 
discrete observation times. 
We define a bi-dimensional Markov chain
by: $X_0 :=(S_0,1)$ and for $1 \leq p \leq P$: $X_p := ( S_{t_p}, (\prod_{q=1}^p
S_{t_q})^{1/p})$, with $S_{t_p}$ given by (\ref{eq:Stp}). Here $s=2$ and $d=1$. 
We choose:   
$S_0=100$,  $K=90$, $r=\log_{10}(1.09)$,  $\sigma=0.2$, $T=240/365$,
 $P= 10$ and $\Delta t_p = T/P$, for $1\leq p\leq P$.  
 
\begin{figure}
\begin{center}
\includegraphics[width=.9\textwidth,height=.4\textheight]{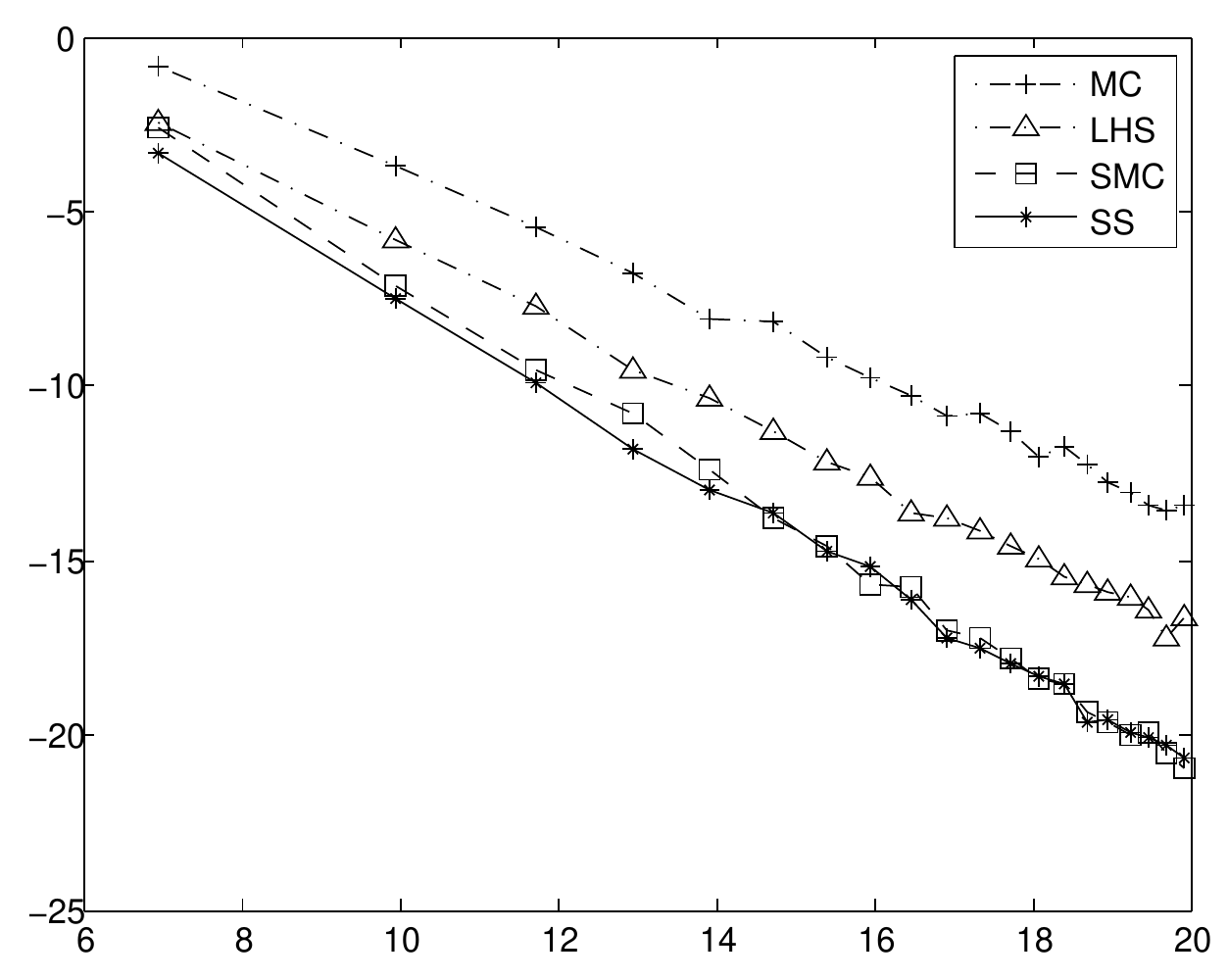} 
\end{center}
\caption{Asian option. Sample variance of $100$ copies of the calculation of $C_{\mathcal{A}}$ as a function of $N$. MC ($+$), LHS ($\triangle$), SMC ($\square$) 
and SS ($*$) outputs, in log-log scale (base $2$).}
\label{fig:asian_var}
\end{figure}

We compute the sample variance of $100$ independent calculations of $C_{\mathcal{A}}$ by MC, LHS, SMC and SS methods. The variances as functions of $N$, for $N=(5m)^3$, $1\leq m\leq20$, are plotted in Figure \ref{fig:asian_var}. 

The order $\beta$ of the variance is estimated using linear regression and the results are given in Table~\ref{tab:asian}.
The convergence rates are not far from those proved for numerical integration in dimension $3$.

\begin{table}
\begin{center}
\begin{tabular}{cccc}
\hline
MC & LHS   & SMC   & SS \\
\hline 
0.99 & 1.12   &  1.40 &  1.33 \\
\hline
\end{tabular}
\end{center}
\caption{Asian option: order $\beta$ of the variance of the calculation of $C_{\mathcal{A}}$.}
\label{tab:asian}
\end{table}

\begin{figure}
\begin{center}
\includegraphics[width=.9\textwidth,height=.4\textheight]{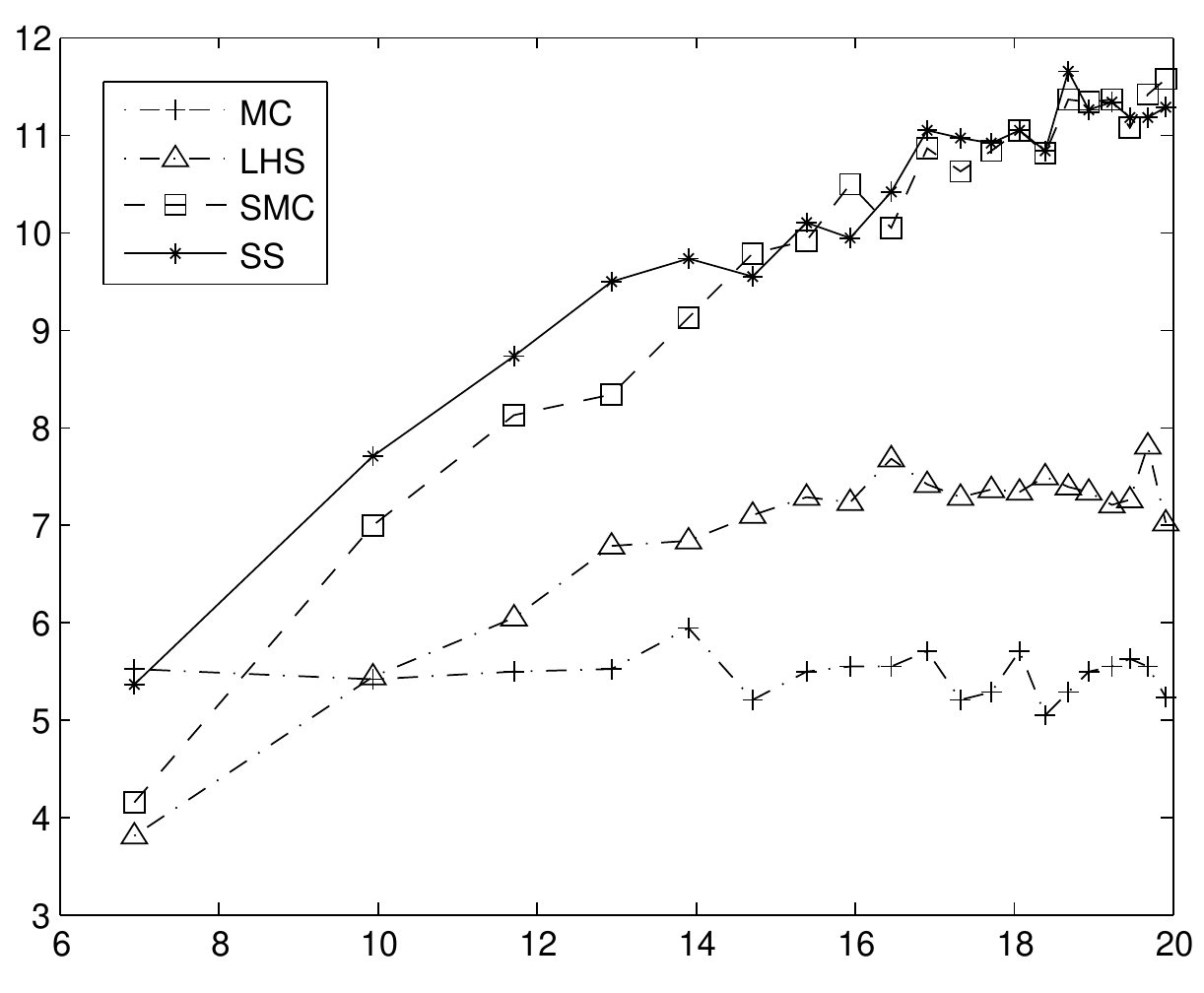} 
\end{center}
\caption{Asian option: efficiency of $100$ copies of the calculation of 
$C_{\mathcal{A}}$ as a function of $N$. Comparison of MC ($+$), LHS ($\triangle$), SMC ($\square$) 
and SS ($*$) outputs, in log-log scale (base $2$).}
\label{fig:asian_eff}
\end{figure}

As before, SMC and SS stratification techniques give smaller variances and 
better convergence rates. But the advantage of the SS algorithm compared to the SMC is lost. At each step of the SS algorithm, the mapping
$\ell \in \{1, \ldots , n \}^3 \to m(w_{\ell}^{\prime}) \in \{1, \ldots , n \} \times \{1, \ldots , n^2 \}$ is not necessarily one-to-one. The efficiencies of the four methods are reported in Figure~\ref{fig:asian_eff}. SMC and SS calculations give similar results and outperform MC and LHS outputs.

\section {Conclusion}
\label{sec:Conclusion}

We have proposed upper bounds for the variance, when we approximate the integral of an indicator function of a subdomain of $I^s$ with stratified Monte Carlo techniques. We have proposed strategies for simulating Markov chains using stratified samples and we have shown on examples that this approach could lead to better efficiency than naive Monte Carlo simulation.

The variance bound of the LHS approximation is obtained in a very restrictive case and should be extended to less specific subdomains of $I^s$. The analysis of stratified simulation of Markov chains remains undone and will be the subject of future work.

\section*{References}


\begin{thebibliography}{[73]}
%

\bibitem {CD89:Problem}
 R. C. H. Cheng, T. Davenport,
 The problem of dimensionality in stratified sampling,  
Management Science  35  (1989) 1278--1296.
%
\bibitem {EFL13:Stratified}
R. El-Haddad,  R. Fakhereddine,  C. L\'ecot,
Stratified {M}onte {C}arlo integration,  in: K. K. Sabelfeld, I. Dimov (Eds.),
Monte Carlo Methods and Applications, De Gruyter, Berlin, 2013, pp. 105--113.
%
\bibitem {EFLV14:Extended}
 R. El-Haddad,  R. Fakhereddine,  C. L\'ecot, G. Venkiteswaran,
Extended {L}atin hypercube sampling for integration and simulation, in: J. Dick, F.Y. Kuo, G.W. Peters, I.H. Sloan (Eds.), Monte Carlo and Quasi-Monte Carlo Methods 2012, Springer, Berlin, 2014, pp.  317--330. 
 %
\bibitem {ELL08:Quasi}
 R. El-Haddad, C. L\'ecot,  P. L'Ecuyer,
Quasi-{M}onte {C}arlo simulation of discrete-time {M}arkov chains on multidimensional state spaces, in: A. Keller , S. Heinrich,  H. Niederreiter (Eds.), Monte Carlo and Quasi-Monte Carlo Methods 2006, Springer, Berlin, 2008,
pp. 413--429.
%
\bibitem{ELLN10:Quasi}
 R. El-Haddad,  C. L\'ecot,  P. L'Ecuyer, N. Nassif,
Quasi-{M}onte {C}arlo methods for {M}arkov chains with continuous multi-dimensional state space, Mathematics and Computers in Simulation 81,
(2010), {560--567}.
%
\bibitem{ES00:Approximating}
M. Evans, T. Swartz, Approximating Integrals via Monte Carlo and Deterministic Methods, Oxford University Press, Oxford, 2000. 
%
\bibitem {Fis96:MonteCarlo} 
G. S. Fishman, Monte Carlo, Springer, New York, 1996.
%
\bibitem{Gla04:MonteCarlo}
P. Glasserman, Monte Carlo Methods in Financial Engineering, Springer,
New York, 2004.
%
\bibitem{Hab66:Modified}
S. Haber, A modified {M}onte-{C}arlo quadrature, 
 Mathematics of Computation 20  (1966) 361--368. 
%
\bibitem{Lec89:Direct}
C. L\'ecot, A {D}irect {S}imulation {M}onte {C}arlo scheme and uniformly distributed sequences for solving the {B}oltzmann equation, Computing 41 
(1989) {41--57}.
%
\bibitem{LC98:Quasi}
{C. L\'ecot, I. Coulibaly},
{A quasi-{M}onte {C}arlo scheme using nets for a linear {B}oltzmann equation},
{SIAM Journal on Numerical Analysis 35} 
(1998)  {51--70}.
%
\bibitem{LT04:Comparison}
C. L\'ecot,  B. Tuffin,
Comparison of quasi-{M}onte {C}arlo-based methods for the simulation
of {M}arkov chains,
Monte Carlo Methods and Applications 10 
(2004)   377--384.
%
\bibitem{LT04:Quasi}
 {C. L\'ecot, B. Tuffin},
 {Quasi-{M}onte {C}arlo methods for estimating transient measures of discrete time {M}arkov chains}, in: H. Niederreiter (Ed.), 
  {Monte Carlo and Quasi-Monte Carlo Methods 2002},
 Springer,  Berlin, 2004, pp. {329--343}.
%
\bibitem{Lec94:Efficiency}
P. L'Ecuyer,
{Efficiency improvement and variance reduction},
  in: J. D. Tew, S. Manivannan, D. A. Sadowski, A. F. Seila (Eds.),
{Proceedings of the 1994 Winter Simulation Conference},
{IEEE Press},
{1994},
 pp.  122--132.
 %
\bibitem{LLL09:Array}
P. L'Ecuyer, C. L\'ecot, A. L'Archev\^eque-Gaudet,
 On array-{R}{Q}{M}{C} for {M}arkov chains: mapping alternatives and convergence rates, in: {P. L'Ecuyer, A. B. Owen (Eds.)},
  {Monte Carlo and Quasi-Monte Carlo Methods 2008},
 Springer, {Berlin},
  {2009},
 pp. {485--500}. 
%
\bibitem  {LLT06:Randomized}
 P. L'Ecuyer, C. L\'ecot, B. Tuffin, 
Randomized quasi-{M}onte {C}arlo simulation of {M}arkov chains with an ordered state space, in: H. Niederreiter,  D. Talay (Eds.)
 Monte Carlo and Quasi-Monte Carlo Methods 2004,
Springer, Berlin, 2006, pp. 331--342.
%
\bibitem{LLT08:Randomized}
P. L'Ecuyer, C. L\'ecot, B. Tuffin,
 A randomized quasi-{M}onte {C}arlo simulation method for {M}arkov chains,
Operations Research 56 
 (2008) 958--975.
%
\bibitem {MBC79:Comparison}
  M. D. McKay,  R. J. Beckman,  W. J. Conover,
A comparison of three methods for selecting values of input variables in the analysis of output from a computer code,
 Technometrics 21 
(1979)  239--245.
%
\bibitem {Nie92:Random}
H. Niederreiter,
Random {N}umber {G}eneration and {Q}uasi-{M}onte {C}arlo {M}ethods,
SIAM, Philadelphia, Pennsylvania, 1992. 
 
%
\bibitem{Owe97:MonteCarlo}
A. B. Owen,
Monte {C}arlo variance of scrambled net quadrature,
SIAM Journal on Numerical Analysis 34 (1997)
1884--1910.
%
\bibitem{Ste87:Large}
M. Stein,
Large sample properties of simulations using {L}atin hypercube sampling,
Technometrics 29  (1987)  143--151.
\end{thebibliography}

\end{document}